\documentclass[amssymb,11pt]{amsart}
\usepackage{amscd,amssymb,euscript,amsthm}
\theoremstyle{plain}
 
\theoremstyle{definition}

\theoremstyle{remark}

\numberwithin{equation}{section}%assigns eqns section numbers
\def\<{\left<}
\def\>{\right>}
\def\cstar{$C^*$-algebra}
\usepackage{pdfsync}
\begin{document}
\title{Helson and subdiagonal operator algebras}
\author{William Arveson}
%
%\thanks{*supported by 
%NSF grant DMS-0100487} 
%
%\address{Department of Mathematics,
%University of California, Berkeley, CA 94720}
%
%\email{arveson@math.berkeley.edu}
%
%\subjclass{46L55, 46L09}
%\subjclass[2000]{Primary 46N50; Secondary 81P68, 94B27}
\date{25 October, 2010}

%\begin{abstract} 
%\end{abstract}

\maketitle

\section{Commutative origins}

Henry Helson is known for his work in harmonic analysis, function theory, invariant subspaces 
and related areas 
of commutative functional analysis.  I don't know the extent to which Henry realized, however, that some of his 
early work inspired significant developments in noncommutative directions -  in the 
area of non-self adjoint operator algebras.  Some of 
the most definitive results were obtained quite recently.  
I think he would have been pleased by that - while vigorously disclaiming any credit.  
But surely credit is due; and in this note I will discuss how his ideas contributed to the noncommutative world of operator algebras.

It was my good fortune to be a graduate student at UCLA in the early 1960s, when 
the place was buzzing with exciting new ideas that had 
grown out of the merger of classical function theory and the more abstract 
theory of commutative Banach algebras as developed by 
Gelfand, Naimark, Raikov, Silov and others.   At the same time, the emerging theory of von Neumann algebras 
and \cstar s was undergoing rapid and exciting development of its own.  One of the directions of that 
noncommutative development - though it 
went unrecognized for many years - was the role of ergodic theory 
in the structure of von Neumann algebras that was pioneered by Henry Dye \cite{DyeGpI}, \cite{DyeGpII}.  
{\em That} Henry would 
become my thesis advisor.  I won't say more about the remarkable 
development of noncommutative ergodic 
theory that is evolving even today since it is peripheral to what I want to 
say here.  I do want to describe the development 
of a class of non-self-adjoint operator algebras that relates to analytic function theory, prediction 
theory and invariant subspaces: Subdiagonal operator algebras.

It is rare to run across a reference to Norbert Wiener's book on prediction theory  \cite{wiePred} in 
the mathematical literature.  That may be partly because the book is directed toward an engineering 
audience, and partly because it was buried as a classified document during the war years.  
Like all of Wiener's books, it is remarkable and fascinating, but not an easy read for students.  
It was inspirational for me, and was 
the source from which I had learned 
the rudiments of prediction 
theory that I brought with me to UCLA as a graduate student.  Wiener was my first mathematical hero.

Dirichlet algebras are a broad class of function algebras that originated in efforts to understand the 
{\em disk algebra} $A\subseteq C(\mathbb T)$ of 
continuous complex-valued functions on the unit circle whose negative Fourier coefficients vanish.  
Several 
paths through harmonic analysis or complex function theory or prediction theory 
lead naturally to this function algebra.  
I remind the reader that a 
{\em Dirichlet algebra} is a unital subalgebra $A\subseteq C(X)$ ($X$ being a compact Hausdorff space) with the 
property that $A + A^* = \{f + \bar g: f,g\in A\}$ is sup-norm-dense in $C(X)$; equivalently, the real parts 
of the functions in $A$ are dense in the space of real valued continuous functions.  
One cannot overestimate the influence of the two papers of Helson and Lowdenslager 
(\cite{hl1}, \cite{hl2})
in abstract function 
theory and especially Dirichlet algebras.  Their main results are beautifully summarized in 
Chapter 4 of Ken Hoffman's book \cite{hofBanSp}.

Along with a given Dirichlet 
algebra $A\subseteq C(X)$, one is frequently presented with a distinguished complex homomorphism 
$$
\phi: A\to \mathbb C  
$$
and because $A+ A^*$ is dense in $C(X)$, one finds that there is a unique probability measure $\mu$ 
on $X$ (of course I really mean unique 
{\em regular Borel} probability measure) that represents $\phi$ in the sense that 
\begin{equation}\label{eq1}
\phi(f)= \int_X f\,d\mu, \qquad f\in A.  
\end{equation}

Here we are more concerned 
with the closely related notion of weak$^*$-Dirichlet algebra $A\subseteq L^\infty(X,\mu)$, in which 
uniform density of  $A+A^*$ in $C(X)$ is weakened to the requirement that 
$A+A^*$ be dense in $L^\infty(X,\mu)$ relative to the weak$^*$-topology of $L^\infty$.  
Of course we continue to require that the linear 
functional (\ref{eq1}) should be multiplicative on $A$.  

\section{going noncommutative}

von Neumann algebras and \cstar s of operators on a Hilbert space $H$ 
are self-adjoint -- closed under the $*$-operation of $\mathcal B(H)$.  But most operator 
algebras do not have that symmetry; and for non-self-adjoint algebras, there was little theory and 
few general principles in the early 1960s beyond the Kadison-Singer paper \cite{ksTriang} on 
triangular operator algebras (Ringrose's work on nest algebras 
was not to appear until several years later).    

While trolling the waters for a thesis topic, I was 
struck by the fact that so much of prediction theory and analytic function theory had been 
captured by 
Helson and Lowdenslager, while at the same time I could see diverse examples of 
operator algebras that seemed to satisfy noncommutative variations of the axioms for weak$^*$-Dirichlet algebras. 
There had to be a way to put it all together in an appropriate noncommutative context
 that would retain the essence of prediction 
theory and contain important examples of operator algebras.  I worked on that idea for a year or two
and produced a Ph.D. thesis in 1964 -- which evolved into a more definitive 
paper \cite{arvAnal}.  At the time I wanted to call 
these algebras {\em triangular};  but 
Kadison and Singer had already taken the term for their algebras  \cite{ksTriang}.     Instead, these later algebras became 
known as {\em subdiagonal} operator algebras.  

Here are the axioms for a (concretely acting) subdiagonal algebra of operators in $\mathcal B(H)$.  It is a 
pair $(\mathcal A, \phi)$ consisting of 
a subalgebra $\mathcal A$ of $\mathcal B(H)$ that contains the identity operator, is closed in the weak$^*$-topology 
of $\mathcal B(H)$, all of which satisfy 
\begin{enumerate}
\item[SD1:]  $\mathcal A+\mathcal A^*$ is weak$^*$-dense in the von Neumann algebra $\mathcal M$ it generates.  
\item[SD2:]  $\phi$ is a conditional expectation,  mapping $\mathcal M$ onto the von Neumann subalgebra 
$\mathcal A\cap\mathcal A^*$.  
\item[SD3:]  $\phi(AB)=\phi(A)\phi(B)$, for all $A,B\in\mathcal A$.  
\end{enumerate}
What [SD2] means is that $\phi$ should be an idempotent linear map from $\mathcal M$ onto $\mathcal A\cap\mathcal A^*$, that carries positive operators to positive operators, is 
continuous with respect to the weak$^*$-topology, and is faithful in the sense that  for every positive operator $X\in \mathcal M$, 
$\phi(X)=0\implies X=0$. 

We also point out that these axioms differ slightly from the original axioms of \cite{arvAnal}, 
but are equivalent when the algebras are weak$^*$-closed.   

Examples of subdiagonal algebras:  
\begin{enumerate}
\item The pair $(\mathcal A,\phi)$, $\mathcal A$ being the algebra of all lower 
triangular $n\times n$ matrices, 
$\mathcal A\cap \mathcal A^*$ is the algebra of diagonal matrices, and 
$\phi: M_n\to \mathcal A\cap \mathcal A^*$ is the map that replaces a matrix with its diagonal part.    
\item Let $G$ be a countable discrete group which can be totally ordered by a relation $\leq$ 
satisfying $a\leq b\implies xa\leq xb$ for all $x\in G$.  There are many such groups, including 
finitely generated free groups (commutative or noncommutative).  Fix such an order $\leq$ on $G$ and 
let $x\mapsto \ell_x$ be the natural (left regular) unitary representation of 
$G$ on its intrinsic Hilbert space $\ell^2(G)$, let $\mathcal M$ be the weak$^*$-closed linear span of all 
operators of the form $\ell_x$, $x\in G$, and let $\mathcal A$ be the weak$^*$-closed linear span of 
operators of the form $\ell_x$, $x\geq e$, 
$e$ denoting the identity element of $G$.  Finally, let $\phi$ be the state of $M$ defined 
by 
$$
\phi(X) = \langle X\xi,\xi\rangle, \qquad X\in M, \quad \xi=\chi_{e}.  
$$
If we view $\phi$ as a conditional expectation from $\mathcal M$ to the 
algebra of scalar 
multiples of the identity operator by way of $X\mapsto \phi(X)\mathbf 1$, then we obtain a subdiagonal 
algebra of operators $(\mathcal A,\phi)$.   
\item There are natural examples of subdiagonal algebras in $II_1$ factors 
$\mathcal M$ that are based 
on ergodic measure preserving transformations that will be familiar to operator algebraists
 (see \cite{arvAnal}).  
\end{enumerate}

In order to formulate the most important connections with function theory and prediction theory, one 
requires an additional property called {\em finiteness} in \cite{arvAnal}:  there should be a distinguished 
tracial state $\tau$ on the von Neumann algebra $\mathcal M$ generated by $\mathcal A$ that preserves 
$\phi$ in the sense that $\tau\circ\phi=\tau$.  Perhaps we should indicate the choice of $\tau$ 
by writing $(\mathcal A, \phi,\tau)$ rather than $(\mathcal A,\phi)$, 
but we shall economize on notation by not doing so.  

Recall that the simplest form of {\em Jensen's inequality} makes the 
following assertion about 
functions $f\neq 0$ in the disk algebra: {\em 
$\log|f|$ is integrable around the unit circle, and the geometric mean of 
$|f|$ satisfies 
\begin{equation}\label{je}
|\frac{1}{2\pi}\int_\mathbb T f(e^{i\theta})\,d\theta| \leq \exp \frac{1}{2\pi}\int_\mathbb T \log|f(e^{i\theta})|\,d\theta. 
\end{equation}
}      

In order to formulate this property for subdiagonal operator algebras we require 
the determinant function of Fuglede and Kadison \cite{fkDet} - defined as follows for invertible 
operators $X$ in $\mathcal M$: 
\begin{equation*}
\Delta(X) = \exp \tau(\log|X|), 
\end{equation*} 
$|X|$ denoting the positive square root of $X^*X$.  
There is a natural way to extend the definition 
of $\Delta$ to arbitrary (noninvertible) operators in $\mathcal M$.  For example, when $\mathcal M$ is the 
algebra of $n\times n$ complex matrices and $\tau$ is the tracial state, 
$\Delta(X)$ turns out to be the positive $n$th root of $|\det X|$.

Corresponding to (\ref{je}), we will say that a finite subdiagonal algebra $(\mathcal A,\phi)$ 
with tracial state $\tau$ satisfies {\em Jensen's inequality} if  
\begin{equation}\label{Je}
\Delta(\phi(A))\leq \Delta(A), \qquad A\in\mathcal A,  
\end{equation}
and we say that $(\mathcal A,\phi)$ 
satisfies {\em Jensen's formula} if 
\begin{equation}\label{Jf}
\Delta(\phi(A)) = \Delta (A),\qquad A\in\mathcal A\cap\mathcal A^{-1}.  
\end{equation}
It is not hard to show that (\ref{Je})$\implies$(\ref{Jf}). 

Finally, the connection with prediction theory is made by reformulating a classical theorem of Szeg\"o, one 
version of which can be stated as follows:  For every positive function $w\in L^1(\mathbb T,d\theta)$ one has  
\begin{equation*}
\inf_f\int_\mathbb T|1+f|^2 w\,d\theta= \exp\int_\mathbb T \log w \,d\theta,   
\end{equation*}
$f$ ranging over trigonometric polynomials of the form 
$a_1e^{i\theta}+\cdots+a_ne^{in\theta}$.
In the noncommutative setting, there is a natural way to extend the definition of determinant 
to weak$^*$-continuous positive linear functionals $\rho$ on $\mathcal M$, and the proper replacement 
for Szeg\"o's theorem turns out to be the following somewhat peculiar statement:  For every weak$^*$-continuous 
state $\rho$ on $\mathcal M$, 
\begin{equation}\label{St}
\inf\rho(|D+A|^2) = \Delta(\rho), 
\end{equation}
the infimum taken over $D\in\mathcal A\cap\mathcal A^*$ and $A\in\mathcal A$ with  
$\phi(A)=0$ and $\Delta(D)\geq 1$.

In the 1960s, there were several important examples for which I could prove properties (\ref{Je}), (\ref{Jf}) and (\ref{St}); 
but I was unable to establish them in general.  The paper \cite{arvAnal} contains 
the results of that effort.  Among other things, it was shown that every subdiagonal algebra 
is contained in a unique {\em maximal} one, and that maximal subdiagonal algebras admit 
{\em factorization}: {\em Every invertible positive operator in $\mathcal M$ has the form 
$X=A^*A$ for some $A\in\mathcal A\cap\mathcal A^{-1}$.}  Factorization was then used to show 
the equivalence of these three properties for arbitrary {\em maximal} subdiagonal algebras.

\section{Resurrection and Resurgence}

I don't have to say precisely what {\em maximality} means because, in an important development 
twenty years later, Ruy Exel \cite{exl} showed that the concept is unnecessary by proving 
the following theorem:  
{\em Every (necessarily weak$^*$-closed) subdiagonal algebra is maximal.}  
Thus, factorization holds {\em in general} and the three properties (\ref{Je}), (\ref{Jf}), (\ref{St}) are 
{\em always} equivalent.  

Encouraging as Exel's result was, 
the theory remained unfinished because no proof existed that Jensen's inequality, for example, 
was true in general.  Twenty more years were to pass before 
the mystery was lifted.  In penetrating 
work of  Louis Labuschagne and David Blecher \cite{labSz}, \cite{bl2}, \cite{bl3},  \cite{bl5} it was shown that, 
not only are the three desired properties true in general, but virtually all of the classical theory 
of weak$^*$-Dirichlet function algebras generalizes appropriately to subdiagonal 
operator algebras.  

I hope I have persuaded the reader that there is an evolutionary path from the original ideas of Helson 
and Lowdenslager, through 40 years of sporadic progress, to a finished and elegant 
theory of noncommutative operator algebras that embodies 
a remarkable blend of complex function theory, prediction theory, and invariant subspaces.

\bibliographystyle{alpha}
\newcommand{\noopsort}[1]{} \newcommand{\printfirst}[2]{#1}
  \newcommand{\singleletter}[1]{#1} \newcommand{\switchargs}[2]{#2#1}

\end{document}